\documentclass[11pt]{article}
\usepackage{graphicx,amsmath,subfigure}
\usepackage{amssymb,float}

\usepackage[default]{jasa_harvard}    
\usepackage{JASA_manu}

\begin{document}

\title{Estimation of relative efficiency of adaptive cluster vs traditional sampling designs applied to arrivals of sharks}
\author{Aneesh S. Hariharan \\
Interdisciplinary Program in Quantitative \\ Ecology and Resource
Management \\ University of Washington, Seattle, WA 98195  \\ 
email: \texttt{aneesh@amath.washington.edu }\\
Vincent F. Gallucci \\
School of Aquatic and Fishery Sciences, Shark Research Lab \\ University of Washington, Seattle, WA 98195  \\ 
email: \texttt{vgallucc@uw.edu} \\
Craig Heberer \\
National Marine Fisheries Service, Southwest Region\\ 501 W. Ocean Blvd., Suite 4200, Long Beach, CA, 90802  \\ 
email: \texttt{craig.heberer@noaa.gov} }

\maketitle

\newpage

\mbox{}
\vspace*{2in}
\begin{center}
\textbf{Author's Footnote:}
\end{center}
Aneesh Hariharan is Doctoral Student,
Interdisciplinary Program in Quantitative Ecology and Resource
Management (QERM).  Mailing address: University of Washington
Loew Hall 304
Box 352182, Seattle, WA 98195 (email: aneesh@amath.washington.edu). Vincent Gallucci is Professor, School of Aquatic and Fishery Sciences and Professor, QERM. Craig Heberer is Fisheries Biologist, National Marine Fisheries Service, Southwest Region, Long Beach, California.

\newpage
\begin{center}
\textbf{Abstract}
\end{center}
Adaptive Cluster Sampling (ACS) is introduced as a technique to use when "natural" groupings are evident in a spatially distributed population, especially sparsely distributed populations. An ACS sampling design will allow efficient allocation of survey manpower with more effective decision rules for where/when to allocate those resources. In particular, given a clustered distribution, ACS would result in a lower variance of the mean estimator than Simple Random Sampling (SRS).  This paper derives a linear inequality based on a ratio of variances and the SRS sample size to determine the conditions under which ACS would be a more appropriate sampling strategy. This inequality could be used to make preliminary decisions on the potential benefits of pursuing an ACS design over an SRS design. The relationship between relative efficiency of ACS over SRS is discussed for (1) variable conditions for neighborhood expansion (Conditions to Adpat) (2) the degree of clustering and (3) the percent of “hits” in a particular neighborhood. Simulation results demonstrate how 'rare' clusters should be for ACS to be a better design than SRS. The relationships that emerge from the simulation provide insights that were not apparent in the analysis. The results are important in the estimation of thresher shark landings along the California coast.

\vspace*{.3in}

\noindent\textsc{Keywords}: {Adaptive Cluster Sampling, Simple random sampling, Realtive efficiency}

\newpage

\section{Introduction}

Natural resource management applications are rich in situations where distributions of items to be sampled occur in patterns that are clustered in space or time. This usually occurs because there are inherent associations (e.g., behavioral, environmental) between items or organisms that may be unknown to the observer or scientist but which essentially determine the properties of that distribution (E.g., means, variances, characteristics of clusters etc). Several sampling designs such as simple random, stratified, systematic, systematic cluster, unequal probability, traditional cluster, adaptive cluster etc. are used to estimate these distribution properties. \cite{mier2008estimating} compared these techniques applied to widespread aggregation of larval walleye pollock in the Gulf of Alaska. In contrast to the preceding, in this paper, we focus on adaptive cluster sampling and attempt to compare it to traditional cluster sampling and simple random sampling. A comprehensive review of adaptive cluster sampling can be found in \cite{christman2000review}, \cite{thompson2004sampling} and \cite{brown2003designing}. \cite{quinn39hoag} developed an adaptive cluster sampling design to estimate the abundance of rockfish populations. Along similar lines, an adaptive sampling design was developed for application to survey sampling of catches of thresher sharks in a Southern California sport fishery \cite{gallucciadaptive}.\\

In the above application, estimators were derived to estimate the abundance of common thresher sharks (Alopias Vulpinus), a U.S. federally managed species (PFMC, 2004), caught in sport fishery primarily off the coast of Southern California. Problems in the estimation of sport fishery characteristics are a different challenge from the estimation of the same properties from commercial fisheries where extensive capture gear is employed. In this particular sport fishery, the California Department of Fish and Wildlife monitors public boat ramps, docks, and marinas where sport vessels return after fishing for thresher sharks. Marinas are located along the coast of California from Point Conception to the Oregon coast. Some of these marinas are clustered close to city areas. Captures are generally from fishing relatively close to the coastline, usually within state water boundaries (i.e., out to three nautical miles). Fishing gear for thresher shark sport fisheries is less of an issue than for  commercial thresher shark fisheries since sport fisheries are  restricted to one fishing rod per individual fisherman, each of whom has a license. There may be multiple anglers and/or rods on any given vessel trip but the Department of Fish and Wildlife samplers query for this. \\

The common thresher shark is a member of the genus Alopidae distinguished from other sharks by the presence of a large asymmetric split tail fin that it uses to stun and capture its preferred prey consisting primarily of mackerel, sardine, and squid. It is classified as a migratory species and thus falls under international laws that regulate migratory species. (Article 64 of the UN Laws of the Seas (UNCLOS)). Therefore, there are national and international responsibilities to monitor the catch of these sharks, bearing in mind that the word ‘catch’ is a euphemism for mortality. The sport fishery targets sharks of all age classes with the predominant catch being juveniles (pre-reproductive) from the first to the fifth year of their lives. During the spring to early summer months, older and actively reproducing thresher sharks will aggregate to feed and when conditions favorable to give birth to pups. There is a demographic relationship between juveniles and the mature adults so that over-harvesting of the young ones may lead to an insufficient number moving into the mature, reproductive category, leading to the possibility of local depletion and over-exploitation , a phenomenon that   occurred back in the late 1980s-early 1990s when the U.S. commercial drift gillnet fishery for swordfish and common thresher sharks was fished at an unsustainable level . Conservation and management measures were put in place to address the sustainability issues leading to a rebound in the population. These conservation measures focused solely on commercial fisheries leaving unchanged the regulations on recreational shark fisheries. These regulations include a 2 shark per day per angler bag limit and no season or size restrictions. \\
Adaptive Cluster Sampling (ACS) is applied to find an estimate of the number of threshers caught in particular time intervals over a range of areas (sampling units consisting of public launch ramps and marinas), in which clustering of boats returning with shark catches have been observed. The sharks have demonstrated a north-south movement pattern along the coast responding to cues such as water temperature, food etc. An earlier study used a stratified simple random sampling design of the public launch ramps and marinas and provided in some instances, estimates with huge variances. In that study, the California coast was split into strata, north-south along the coast; within each stratum random sites were chosen for sampling. The ratio: the number of randomly chosen sites per stratum/ the total number of sites in that stratum (sampling proportion per stratum) across all strata was held constant, to derive the necessary estimators.\\

Traditional cluster sampling is an alternative design to stratified simple random sampling to estimate catch. In this method, the estimate for the variance of the mean number of boats returning with shark catch will require prior knowledge of the number of public ramp and marina clusters and consequently the number of public ramps and marinas in each cluster, which would then be subsampled to get the total number of boats with shark catches. However, criteria to define clusters are not available, thus sample designs are not necessarily unique. 
The characteristics of ACS that fit the circumstances not achieved by cluster sampling are (a) it adapts to whether selected sample units fall within a cluster and (b) adapts to detecting a cluster, if one exists, without apriori knowledge of their existence . Achieving a theoretical comparison of the efficiency of traditional cluster sampling and ACS seems impractical, because of the fundamental difference in implementing the two designs. In ACS, one does not know apriori the sample size needed to achieve a certain level of variance, since the sample size evolves as the sampling activity progresses. In contrast, traditional cluster sampling needs complete knowledge about the number of clusters and within each cluster, the primary sampling units, to achieve a certain a level of the variance of an estimate. \\

Appendix A presents estimators from traditional cluster sampling and clarifies why it is impractical to make a theoretical comparison between traditional cluster sampling and ACS.  ACS is equivalent to SRS if no clusters are found. Thus, comparison of ACS to SRS is straightforward and more natural in terms of the mathematics while comparison of ACS to traditional cluster sampling would be natural with respect to clusters. In this paper, the authors compare ACS vs SRS and include a simulation study that provides insights into the sensitivity and efficiency of ACS relative to the degree of clustering.\\

The outstanding advantage of adaptive sampling is its suitability for sparsely distributed, clustered populations. When comparing sampling designs, the efficiency of one design over another is measured in terms of relative variances. Another important component of efficiency is measured in terms of the cost to carry out the sampling design. The relative efficiency addressed in this paper is in terms of variances only, which is the basis for comparing ACS to SRS sampling designs.

\section{METHODOLOGY }
\label{sec:Adm}

In an adaptive cluster sampling design, an initial random sample of $n_1$ sample units is selected from $N$ units of the population. Each unit in the population may represent a quadrat or some fixed region with an associated neighborhood. The neighborhood of sample unit $i$  is defined as unit $i$ plus all the adjacent units that satisfy a certain threshold criterion to determine if additional samples should be taken in that neighborhood. Whenever the $y$-value of unit $i$ in the sample satisfies a criterion  (e.g. $y_i>C$ ), where $C$ is a specified constant), all units in the neighborhood of $i$ are added to the sample. This process is iterated for every sample unit that satisfies the threshold criterion, including the ones that are newly added to a network, where a network is defined as a subset of all sample units within a cluster such that if any sample unit of the network is selected, all other plots of this network will enter into the sample. For rare event sampling, such as the theresher shark application, a criterion $y_i>1$ is used for expanding neighborhoods, since the problem involved sampling for a rare species, thus even a small criterion was enough 
to trigger an adaptive sampling network. As can be clearly seen, specification of a criterion will vary from one study to another and is usually done in consultation with people familiar with the system and the response variable and its behavior. It is noted here that just as stratified SRS improves precision over SRS, similarly ACS with stratification can reduce the variance over a standalone ACS design.

\section{DEFINITION OF SYMBOLS USED}
The symbols used in this paper are in concurrence with the symbols used in \cite{thompson1996adaptive} and \cite{thompson1992sampling}.\\
Let $N$ - Number of units in a finite population with a variable of interest $y_i$ in the  $i^{th}$ unit.\\
$n_1$ - Number of initial random sample units chosen from  units in the population for the Adaptive Cluster Sampling (ACS) Design.\\
$K$ - Number of networks in the population.\\
$y_j$ - Number of observations of the variable of interest in sample unit $j$ .\\
$k(i)$ - Label for the network that includes unit $i$.\\
$B_k$ - Set of units that are present in the $k^{th}$ network.\\
$m_{k(i)}$ - Number of units in network $B_{k(i)}$.\\
$w_{k(i)}$- Average of the $y$- values of the units in the network that include unit $i$ . (It is easily seen that $w_{k(i)}$ =$\frac{1}{m_{k(i)}}\Sigma_{j \in B_k}y_j$ ).\\
$m$ - Number of sample units for a simple random sample where units are chosen randomly from N units.\\
$\overline{var(acs)}$ – Variance estimator for the estimated mean as a result of adaptive sampling.\\
$\overline{var(srs)}$- Variance estimator for the estimated mean as a result of simple random sampling.

\section{Derivations}
Thompson (1996, pp. 151, eq. 5.1)  provides an expression for the relative efficiency of a conventional simple random sampling (SRS) design with sample size $m$, to the adaptive cluster sampling  (ACS) strategy of initial sample size $n_1$ in terms of their respective variances. This is given by:
\begin{equation}
\frac{\overline{var(acs)}}{\overline{var(srs)}}=\frac{m}{n_1}\frac{N-n_1}{N-m}\big[1-\frac{\Sigma_{k=1}^{K}\Sigma_{i \in B_k}(y_i-w_{k(i)})^2}{\Sigma_{i=1}^{N}(y_i-\mu)^2}\big]
\end{equation}
Thus, adaptive sampling will have a lower variance than simple random sampling if the within-network variance of the population is sufficiently high and occurs in patches.  Moreover, for given sample sizes m and $n_1$, adaptive cluster sampling will be more efficient if the within network variation represents a large proportion of the overall variation, since larger the proportion,  smaller the quantity within equation (1). The adaptive criterion $C$  that defines the conditions for expansion of a neighborhood plays a critical role here. This criterion governs the partitioning of networks in the population and is directly proportional to the characteristics of highly aggregated clusters that occur in the population. 
Equation (2) below describes the conditions under which ACS will be a "superior" sampling design compared to SRS, ie ACS will lead to a lower variance than SRS. The derivation for this inequality of provided in Appendix B.
\begin{equation}
\big(\frac{1}{n_1}-\frac{1}{m}\big)\sigma^2<\frac{N-n_1}{n_1N}\frac{\Sigma_{k=1}^{K}\Sigma_{i \in B_k}(y_i-w_{k(i)})^2}{N-1}
\end{equation}
Immediate observations from Equation (2):\\
\begin{enumerate}
\item If $n_1=m$, then the variance of estimators obtained using adaptive sampling is always lower than the variance obtained using simple random sampling.  This is because the left-hand side (LHS) becomes 0 and the right-hand side (RHS) is always a strictly positive quantity.
\item If $n_1>m$, then the LHS is a negative number and the RHS is a positive number. Thus, adaptive sampling will perform better than SRS.
\item Adaptive sampling will require a smaller initial sample to obtain the same precision as SRS. That is, $n_1<m$  .
\end{enumerate}
In general, for a complex survey design, we usually need a larger sample size to generate the same level of precision as SRS. It is only the initial number of samples for adaptive sampling that is claimed to be lower than those required for SRS. At the end of the study, with expanded neighborhoods grouped in clusters, the number of sample units sampled is higher. Thus, we may achieve better precision using adaptive cluster sampling than simple random sampling for a fixed cost, since traveling from one site to another is not random and depends upon the neighborhood of the site visited on the previous sample day. That is, even if a larger number of sample units are traversed, the cost may be lower for the same, or greater level of precision as SRS.
Equation (A.7) in Appendix B expresses the variance of the population from two different sampling designs, the one from SRS and the double summation proportional to the ACS variance. The relationship between these variances and the variance of the means follows. For SRS, it is $\frac{\sigma^2}{m}$ , and for ACS, it is given by equation(4) where $\tilde{\mu}$ is the mean estimator for ACS.
The sum of squares for the response variable can be considered to be a sum of within network and between network components, similar to a one way ANOVA (Thompson 1996) i.e
\begin{equation}
\Sigma_{i=1}^{N}(y_i-\mu)^2=\Sigma_{k=1}^{K}\Sigma_{i \in B_k}(y_i-w_{k(i)})^2+\Sigma_{i=1}^{N}(w_{k(i)}-\mu)^2\noindent
\end{equation}
Therefore the variance of the adaptive cluster sampling estimator   is found from the above by multiplying both sides of the equation by $\frac{N-n_1}{n_1N(N-1)}$  and rearranging terms, i.e.
\begin{equation}
Var(\tilde{\mu})=\frac{N-n_1}{n_1N(N-1)}[\Sigma_{i=1}^{N}(y_i-\mu)^2-\Sigma_{k=1}^{K}\Sigma_{i \in B_k}(y_i-w_{k(i)})^2]\noindent
\end{equation}
Rearranging equation(4) leads to 
\begin{equation}
\frac{N-n_1}{n_1N(N-1)}\Sigma_{k=1}^{K}\Sigma_{i \in B_k}(y_i-w_{k(i)})^2=\frac{N-n_1}{n_1N(N-1)}\Sigma_{i=1}^{N}(y_i-\mu)^2-Var(\tilde{\mu})
\end{equation}
Plugging the right hand side of the above into (2), we get
\begin{equation}
(\frac{1}{n_1}-\frac{1}{m})\sigma^2<\frac{N-n_1}{n_1N(N-1)}\Sigma_{i=1}^{N}(y_i-\mu)^2-Var(\tilde{\mu})
\end{equation}
i.e. substituting the srs population variance,$\sigma^2$  
\begin{equation}
(\frac{1}{n_1}-\frac{1}{m})\sigma^2<\big(\frac{N-n_1}{n_1N}\big)\sigma^2-Var(\tilde{\mu})
\end{equation} 
Rearranging and factoring, we get
\begin{equation}
\big(\frac{1}{m}-\frac{1}{N}\big)\sigma^2>Var(\tilde{\mu})
\end{equation}
which can be rearranged to highlight the ratio of variances as:
\begin{equation}
\frac{1}{m}>\frac{1}{N}+\kappa
\end{equation}
where $\kappa=\frac{Var(\tilde{\mu})}{\sigma^2}$\\
If $N<m$, the inequality cannot be valid. In the event of total enumeration, i.e., when $N=m$, variances are zero, which is the limiting case of the inequality. Thus, irrespective of the initial random sample size chosen for adaptive sampling, if (8) and (9) hold, ACS performs better than SRS. The last option of $m<N$  is more complex and can be expressed in terms of either the $N,m$ relationship or the ratio of variances.
To further investigate the   relationship in (9), multiply both sides of (9) by $m$. Since $\kappa=\frac{Var(\tilde{\mu})}{\sigma^2}$  and since $\frac{\sigma^2}{m}$  is the variance of the sample mean from SRS (without the finite population correction (FPC)), (9) can be rewritten as 
\begin{equation}
1>\frac{m}{N}+\frac{Var(\tilde{\mu})}{(\frac{\sigma^2}{m})}=\frac{m}{N}+\frac{Var(\tilde{\mu})}{Var(\bar{y_m})}
\end{equation}
i.e.
\begin{equation}
1>\frac{m}{N}+\kappa_1
\end{equation}
where $\kappa_1=\frac{Var(\tilde{\mu})}{Var(\bar{y_m})}=m\kappa$\\
Equation (10) can be rearranged as:
\begin{equation}
N(1-\kappa_1)>m
\end{equation} 
The “linear inequality” (12) partitions the space defined by the SRS sample size ($m$) and the ratio of variances of the means of ACS and SRS,$\kappa_1$. The line (when (12) is an equality) has a vertical intercept of population size ($N$) (total enumeration $m=N$) and a negative slope of $–N$. This defines two regions, $N-N\kappa_1=m$ in a ($\kappa_1,m$) space. The region above this line $m>N$ is impossible. $\kappa_1=0$ implies $m=N$ (total enumeration) and, where $m<N$ with $\kappa_1<1$ implies the feasible region. This region, below the line defines the feasible solution set of all combinations of m and the ratio mean variances, where ACS is a superior sampling strategy. 

\section{Simulation of sampling on spatial data}
The theory of spatial distributions makes use of an index, the Variance to Mean Ratio (VMR) to describe different spatial patterns. VMR indices are also known as:  index of dispersion, dispersion index, coefficient of dispersion, and the coefficient of variation. The VMR is a normalized measure of the dispersion of realizations of a probability distribution. For example, standard statistical models such as Poisson, Binomial etc. are often associated with particular random patterns. Mathematically, the VMR is defined as the ratio of the variance $\sigma^2$ to the mean $\mu$ of the random variable associated with a spatial distribution.
The Poisson distribution has equal variance and mean, giving it a $VMR = 1$. The geometric and negative binomial distributions each have a $VMR > 1$, while the binomial distribution has $VMR < 1$, and the uniform distribution has $VMR = 0$. Each of these distributions’ VMR is associated with a spatial pattern described in terms of the “dispersion” of the realization. This relationship between VMR and the spatial distributions is summarized and described as follows:
\begin{table}[h!]
\caption{Relationship between distributions, VMR and their respective descriptions}
\begin{center}
  \begin{tabular}{| c | c | c |}
    \hline
    Distribution & VMR & Description\\ \hline
    Uniform & 0 & not dispersed  \\ \hline
    Binomial & $0<VMR<1$ & under dispersed\\ \hline
    Poisson & VMR=1 & random\\ \hline
    Negative binomial & $VMR>1$ & over dispersed\\ \hline
    \end{tabular}
\end{center}
\end{table}

Spatial patterns can be generated from simulations based on the negative binomial and Poisson pdfs, since they would generate over dispersed data which in turn correspond to clustering.  The aim is to study relative precision (defined as the ratio of $\frac{Var(SRS)}{Var(ACS)}$ as data become more clustered. In our application, let a “hit” be the return of a vessel to a particular marina with at least one shark capture on a particular day. Different spatial distributions of “hits” can be simulated to correspond to spatial patterns of different marinas with hits. For completely unclustered data, the Variance/Mean ratio is 1 (Poisson). As this ratio is increased (data simulated from a negative binomial), the simulated data begin to show clustering. The spatial distribution of recreational fishing boats in the marinas can assume different forms, which, in this paper, range from random to clustered. This could represent, e.g., the marinas distributed along the California coastline, where anglers return with thresher shark catch. \\
Once data is simulated, both SRS and ACS can be applied to these data. A similar study was conducted by \cite{ojiambo2010efficiency} et. al. to estimate plant disease incidence for varying Conditions to Adapt (CA), where CA is a preset criterion (number of hits) to be satisfied before further samples are taken in the neighborhood. If the numbers of hits are below the CA, then the next days' neighborhood samples are determined by SRS. The conditions to adapt are typical of a rare event, e.g. number of thresher shark catches on a particular sample day. It has been shown that ACS is most effective for rare events where the response variable of interest (here, number of “hits”) is infrequently observed (Thompson 1996; Pg. 8,9). Therefore, as data become more clustered (indicated by increasing V/M ratio), the relative precision can be shown to be higher for a lower condition to adapt (e.g., CA=1). It can also be seen that as the “hit level” increases, relative precision decreases (i.e, if more boats return with a thresher shark catch per day, the relative precision decreases). This implies that relative precision is higher for rare, clustered data with a low condition to adapt.
To summarize, the results as addressed in \cite{ojiambo2010efficiency}, which is also applicable in this shark application are :
\begin{enumerate}
\item Relative precision increases as data become more clustered (indicated by an increasing V/M ratio).
\item	Relative precision is inversely related to Conditions to Adapt (CA), i.e., the lower the CA, the higher the relative precision (e.g., if CA=1, the relative precision is the highest). 
\item Relative precision decreases as the “hit level” increases. That is, if more boats return with a thresher catch per day, the relative precision decreases.
\end{enumerate}
\section{ACS vs SRS for varying spread of clusters}
Adaptive cluster sampling (ACS) has been recommended in situations where sparsely distributed populations occur in clusters or special groupings.  Intuitive examples generally include animal populations found in pods (orcas), prides (lions), packs (wolves), and the thresher shark instance that has been discussed as part of this paper.  In each of these animal examples, small groups of animals are likely to avoid other such groups, maintaining large distances between groups. \\
One can imagine an extreme example, where the centers of two or more clusters are so close relative to their spread that the overall population seems more evenly or randomly distributed.  In such an instance the clusters would almost entirely overlap.  For such an extreme example, it would seem that simple random sampling (SRS) would be the preferred design; although ACS may still yield slightly more precise estimates, the additional cost and complexity incurred by ACS may deter use of this scheme.This section begins an examination into whether there is a gradient of relative performance between ACS and SRS given varying relative distances between clustered subpopulations.  Data were simulated in R and both sampling schemes executed; population and variance estimates for both methods are compared.\\
A two-dimensional sampling frame of $N=400$ units (20 by 20) was set up, and x and y coordinates of 5 cluster centers were randomly chosen using a uniform distribution.  In order to simulate varying spread of clusters, a bivariate normal distribution was used to generate 50 data points centered around each of the 5 cluster centers.  Five different spreads of clusters were analyzed: the bivariate normal standard deviation was set at $\frac{2}{3}$ unit, 1 unit, $\frac{3}{2}$ units, 2 units, and 3 units respectively.  This scheme resulted in between 200 and 250 data points within the frame; as some of the 250 points sometimes fell outside of the sampling frame, these were not considered.  Figure 1 shows the data points generated for each of the 5 standard deviations used.\\
For each of the 5 simulated populations, n=10 initial random samples were chosen 100 times.  For the simple random sampling scheme, counts were taken for these 10 samples and estimates calculated.  For the adaptive cluster sampling scheme, the 10 initial samples were used to build an ACS sample, where $C>0$ in a unit triggered sampling of the adjacent 4 units, as outlined in \cite{ojiambo2010efficiency}. Calculated estimates for ACS were compared to those found using SRS for each of the 100 sampling experiments and are averaged and presented in Table 2. The R-code for Table 2 and Figure 1 are in Appendix 3.\\
Examination of Figure 1 shows marked differences between clusters simulated using different variance parameters from a bivariate normal distribution. The first two figures in Figure 1 (SD=2/3 and SD=1) look like what we would expect from clustered animal populations, the last two (SD=2 and SD=3) seem more randomly distributed than clustered and the middle figure (SD=3/2) is somewhere in between.\\ 
Comparing the relative efficiency estimates in Table 2, the improvement in ACS’s performance relative to SRS does improve with tighter, more widely spaced clusters as is to be expected (approximately twice for the most clustered populations).  Variance estimates over 100 sampling experiments for the two least clustered populations were only ~40-50 percent better using ACS than SRS; given the additional cost of ACS it might be more effective in these examples to use SRS. The ‘in between’ population, where the standard deviation of the bivariate normal distribution used to generate the data was 1.5 sampling units, shows that ACS might be a better design since it performs 75 percent better. It is also noted that these results are the result of one replicated simulated study and replication of the same study sometimes yielded only 20 percent improvement of ACS over SRS for the least clustered data. The question remains, is such a reduction worth the extra cost and complexity of adapting ACS?  While this paper does not explore a cost effective approach, it is well worth additional exploration and should be one of the more important questions facing individual study designs.
\begin{table}[H]
\caption{ACS and SRS estimates for various simulated populations. Relative precision increases as clusters become “rare” and ACS is a better design. Estimates for 10 initial samples, repeated 100 times}
\includegraphics[width=180mm,height=220mm]{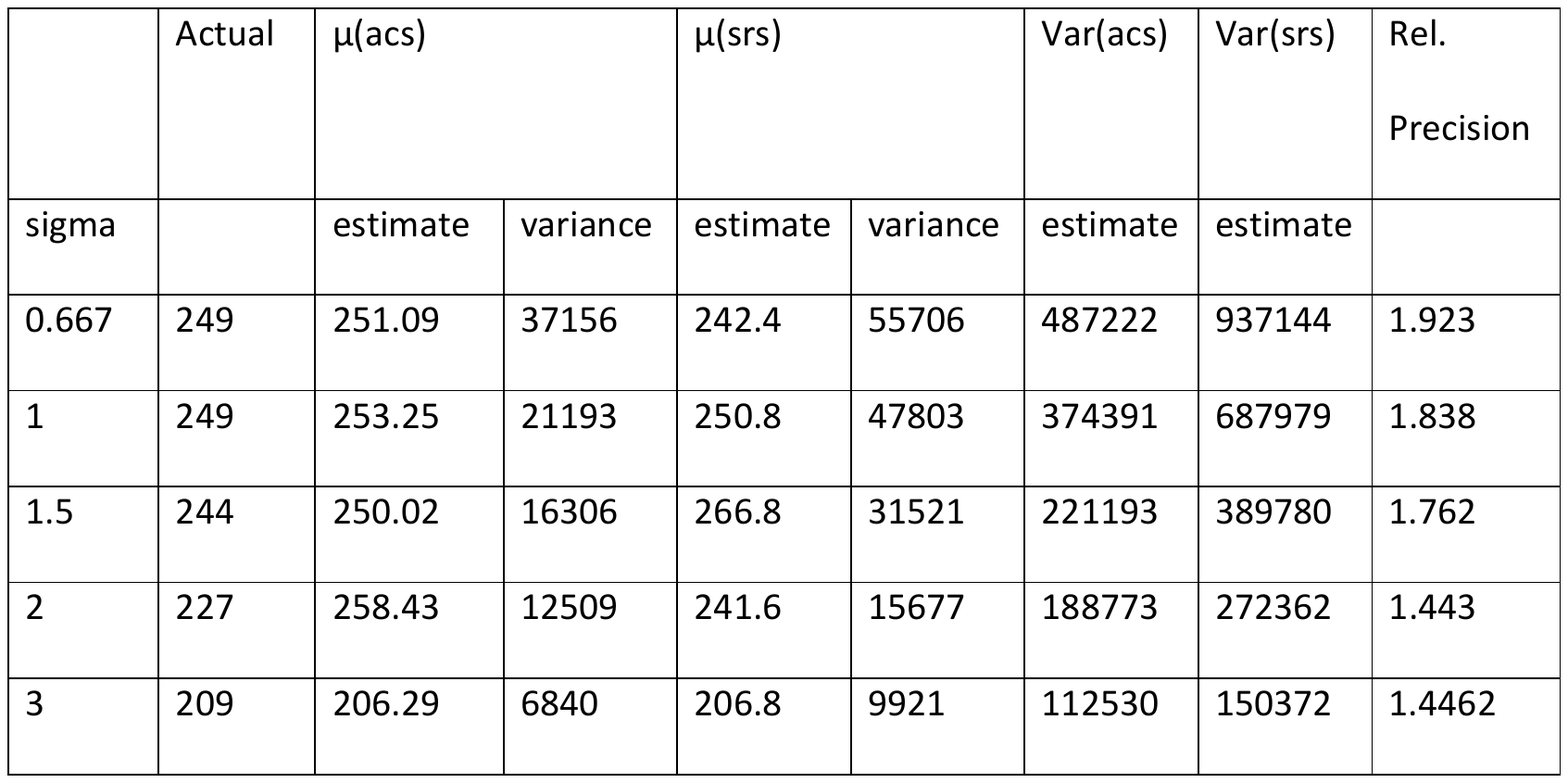}
\end{table}
\nopagebreak
\begin{figure}[H]
\begin{center}
\centering
\subfigure{
\includegraphics[scale=0.25]{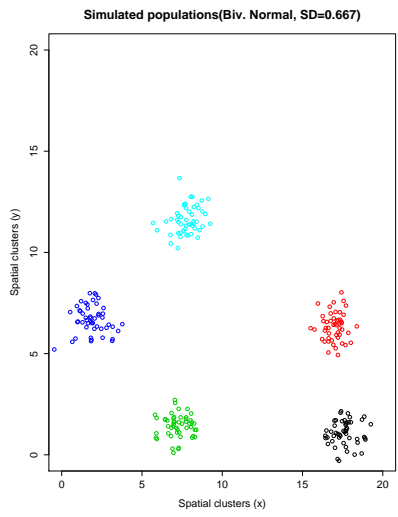}}
\hfill
\subfigure{
\includegraphics[scale=0.25]{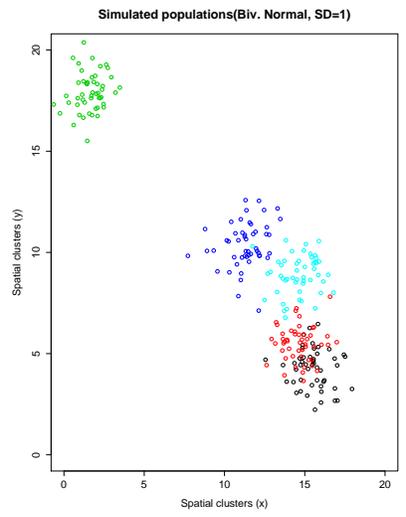}}
\hfill
\subfigure{
\includegraphics[scale=0.25]{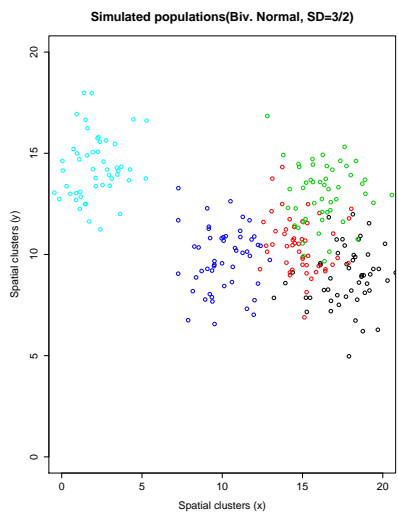}}
\hfill
\subfigure{
\includegraphics[scale=0.25]{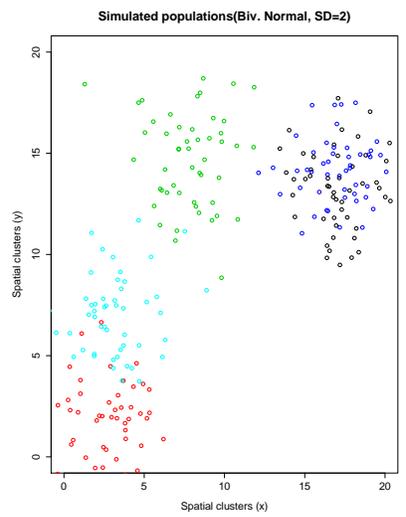}}
\hfill
\subfigure{
\includegraphics[scale=0.25]{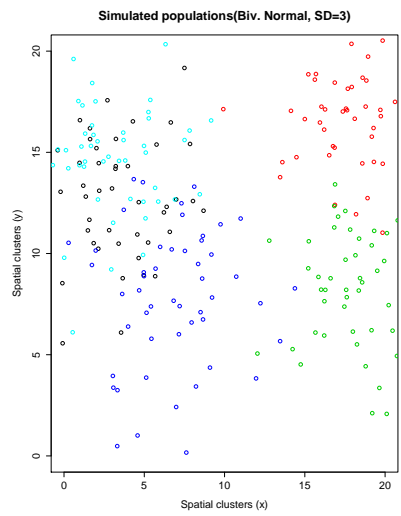}}
\hfill
\caption{Spatial distribution of simulated populations for varying degrees of clustering}
\end{center}
\end{figure}
\section{Discussion}
Adaptive Cluster Sampling differs from traditional sampling designs because in ACS it is impossible to apriori pick an initial sample size ($n_1$) to achieve a preset variance of sampling estimators. Thus, fixing $n_1$ as seen in inequality (2) is not practical. However, use of the ratio of the variance estimators for each sampling design is useful, since it allows one to preset the efficiency of ACS relative to SRS, in terms of:  $n_1$, the SRS sample size and, the desired precision (variance) ratio, $\kappa_1$ (Eqn. 12).\\
Recalling the well-known inverse relationship between variance and sample size, which translates to sampling cost, an adaptive design is more complicated than conventional designs and probably more costly to execute.  ACS will perform better (lower variance estimates relative to SRS) for patchily distributed populations (Thompson, 1996) and, would perform less well in a population where the sample units are not clustered. That is, where sample units are either randomly or uniformly distributed.\\
Having addressed ACS relative to variance and sampling design, consideration of results relative to the application to sampling for rare events, is useful. How rare is rare has been discussed by simulating populations and noting that for tighter clusters, ACS reveals almost twice the efficiency of a traditional SRS design. For unclustered populations or when clusters significantly overlap, ACS still proves to be a better design, but given the additional cost, it involves further exploration on its relative merits.\\
Adaptive sampling designs are, in fact, not found in many sampling books and courses, although, in at least one case, a whole book is devoted exclusively to adaptive sampling designs \cite{thompson1996adaptive}. Nevertheless, in teaching adaptive sampling in other contexts, the second author has often had to provide guidance as to when the selection of adaptive cluster sampling designs will improve over cluster or simple random sampling designs. Referring to equation (2) was the best that could be done. The authors believe the result obtained in Equation (12) and its ramifications is a step forward. Further insight follows from the simulation results.
\section{Acknowledgment}
AH wishes to acknowledge John Skalski, whose sampling theory course at the School of Aquatic and Fishery Sciences and Center for Quantitative Science University of Washington, helped foster ideas presented in this paper. AH also wishes to acknowledge Adam Pope for valuble inputs in the simulation section. Both authors acknowledge the assistance of Steven Stohs, Kirk Lynn, Leanne Laughlin, Suzanne Kohin and others, who presented a thresher shark catch estimation application of adaptive sampling to them and who patiently waited while the necessary theory for this application of adaptive cluster sampling was developed.
\bibliography{example}

\makeatletter   
 \renewcommand{\@seccntformat}[1]{APPENDIX~{\csname the#1\endcsname}.\hspace*{1em}}
 \makeatother

\bibliographystyle{ECA_jasa}
\appendix
\section{Variance estimator for the mean in traditional cluster sampling}
ACS may be confused with traditional cluster sampling, but they are quite different.  They are similar is that both refer to a spatial configuration consisting of clusters of sample units which contain the response variable of interest. In this paper, as well as the primary source (Thompson 1996), ACS is not typically compared to traditional cluster sampling, but rather to SRS. And in typical sampling textbooks, traditional cluster designs are also compared to SRS designs, Thus SRS becomes the baseline design for comparison in Section 4. This is legitimate because most texts compare traditional cluster sampling to SRS. This appendix essentially motivates why traditional cluster sampling cannot be directly compared to ACS.
The following symbols are used in the expressions involving the means and variance estimators of traditional cluster sampling.\\
Let,\\
$N$-Number of clusters in the population.\\
$M_i$-Number of sample units in cluster $i$, $i=1,2,....N$.\\
$Y_{ij}$- value of the response variable $Y$ in sample units $j$ of cluster $i$.\\
\subsection{Estimates at the level of the individual clusters $M_i$}
$M_0=\Sigma_{i=1}^{N}M_i$-Total number of sample units in the population.\\
$\bar{M}=\frac{M_0}{N}$-Mean number of sample units per cluster.\\
$Y_j=\Sigma_{j=1}^{M_i}Y_{ij}$-Total value of the response variable $Y$ in cluster $i$.\\
$\bar{Y_i}=\frac{Y_i}{M_i}$-Mean value of the response variable $Y$ in the sample units of cluster $i$.\\
\subsection{Population parameters at the level of all clusters in the population $N$}
$Y=\Sigma_{i=1}^{N}Y_i$-Total value of the response variable $Y$ of all the sample units in the population.\\
$\bar{Y}=\frac{Y}{N}$-Mean value of response variable $Y$ per cluster.\\
$\bar{\bar{Y}}=\frac{Y}{M_0}=\frac{Y}{N\bar{M}}$-Mean value of response variable $Y$ per sample unit.
\subsection{Sampling estimators from a cluster design}
The mean of the response variable per sample unit is:\\
$\hat{\bar{\bar{Y}}}=\frac{\Sigma_{i=1}^{n}y_i}{n}$, which has sampling variance given by:\\
$V[\hat{\bar{\bar{Y}}}]=\frac{V[\hat{Y}]}{M_0^2}$ where $V[\hat{Y}]=\frac{M_0^2}{n(n-1)}\Sigma_{i=1}^{n}(\bar{y_i}-\bar{\bar{y}})^2$\\
From the expressions above, it is clear that to derive estimators (means and variances) for traditional cluster sampling, one needs complete knowledge of the number of clusters, the number of clusters sampled and the number of primary sampling units per cluster. Thus, traditional cluster sampling can be compared with SRS since, from a design perspective, they are analogous to each other.\\
In contrast, ACS does not allow the estimation of the sample size needed to achieve a variance of the estimate of the mean, rather the sample size evolves from the sampling activity. This paper attempts to circumvent this issue by deriving equation (12), that allows to compare the efficiency of ACS vs SRS, by getting rid of the final ACS sample size from the equation.
One consequence of the above relationship that has not been addressed in this paper, is that the cost of the sampling designs are different. In ACS, cost depends on the number of sample units sampled in a network, which is difficult to know apriori. In some cases, the ACS design will be logistically less expensive since successive samples in a neighborhood may be spatially close to each other, compared to a samples selected by SRS, where the sample units maybe spatially distant from each other. This paper is not meant to address this issue of cost, and focuses on the efficiency of sampling designs purely from the perspective of relative variances.

\section{Inequality to determine when ACS will outperform SRS}
If adaptive cluster sampling outperforms SRS, then
\begin{equation}
\frac{\overline{Var(ACS)}}{\overline{Var(SRS)}}<1
\end{equation}
This implies from Equation 1, that:
\begin{equation}
\frac{m}{n_1}\frac{N-n_1}{N-m}\big[1-\frac{\Sigma_{k=1}^{K}\Sigma_{i \in B_k}(y_i-w_{k(i)})^2}{\Sigma_{i=1}^{N}(y_i-\mu)^2}\big]<1
\end{equation}
i.e.
\begin{equation}
[1-\frac{\Sigma_{k=1}^{K}\Sigma_{i \in B_k}(y_i-w_{k(i)})^2}{\Sigma_{i=1}^{N}(y_i-\mu)^2}\big]<\frac{n_1}{m}\frac{N-m}{N-n_1}
\end{equation}
i.e.
\begin{equation}
1-\frac{n_1}{m}\frac{N-m}{N-n_1}<\frac{\Sigma_{k=1}^{K}\Sigma_{i \in B_k}(y_i-w_{k(i)})^2}{\Sigma_{i=1}^{N}(y_i-\mu)^2}
\end{equation}
i.e.
\begin{equation}
\frac{mN-n_1N}{m(N-n_1)}<\frac{\Sigma_{k=1}^{K}\Sigma_{i \in B_k}(y_i-w_{k(i)})^2}{\Sigma_{i=1}^{N}(y_i-\mu)^2}
\end{equation}
Dividing both sides by $n_1N(N-1)$ and cross multiplying, we get
\begin{equation}
\big(\frac{m-n_1}{mn_1}\big)\frac{\Sigma_{i=1}^{N}(y_i-\mu)^2}{N-1}<\big(\frac{N-n_1}{n_1N}\big)\frac{\Sigma_{k=1}^{K}\Sigma_{i \in B_k}(y_i-w_{k(i)})^2}{N-1}
\end{equation}
Note that the terms not within the parantheses on the LHS represent the SRS population variance with mean $\mu$. Thus,
\begin{equation}
\big(\frac{1}{n_1}-\frac{1}{m}\big)\sigma^2<\big(\frac{N-n_1}{n_1N}\big)\frac{\Sigma_{k=1}^{K}\Sigma_{i \in B_k}(y_i-w_{k(i)})^2}{N-1}
\end{equation}
Here $\sigma^2$ is the population variance. The equation is similar to \cite{thompson1992sampling}pp. 275.
\section{R code for simulating clustered populations}
\begin{verbatim}
library(fMultivar)
xmax<-20
ymax<-20
N<-xmax*ymax
centers<-cbind(round(runif(5,1,xmax-1),2),round(runif(5,1,ymax-1),2))
cluster.1<-rnorm2d(50)*1+c(rep(centers[1,1],50),rep(centers[1,2],50))
cluster.2<-rnorm2d(50)*1+c(rep(centers[2,1],50),rep(centers[2,2],50))
cluster.3<-rnorm2d(50)*1+c(rep(centers[3,1],50),rep(centers[3,2],50))
cluster.4<-rnorm2d(50)*1+c(rep(centers[4,1],50),rep(centers[4,2],50))
cluster.5<-rnorm2d(50)*1+c(rep(centers[5,1],50),rep(centers[5,2],50))
simul.data.2<-rbind(cluster.1,cluster.2,cluster.3,cluster.4,cluster.5)
n<-10
frame<-data.frame(rep(0,N))
frame$x<-trunc((row(frame)-1)/20)
frame$y<-row(frame$x)-20*frame$x-1
for (i in 1:length(frame[,1])) 
{frame[i,1]<-length(which(simul.data.2[,1]>=frame$x[i] & 
simul.data.2[,1]<frame$x[i]+1 & 
simul.data.2[,2]>=frame$y[i] &
simul.data.2[,2]<frame$y[i]+1))                    
}
sum(frame[,1])
mu.acs<-rep(0,100)
var.acs<-rep(0,100)
mu.srs<-rep(0,100)
var.srs<-rep(0,100)

for (h in 1:100) {
  rand.samp<-sample(N,n)
  samp.frame<-matrix(1:N,xmax,ymax)
  rand.ind<-arrayInd(which(samp.frame %in% rand.samp), dim(samp.frame))
  colnames(rand.ind)<-c("x","y")
  names(frame)<-c("count","x","y")
  srs.2<-frame$count[rand.samp]
  count<-matrix(0,xmax,ymax)
for (i in 1:xmax) {
  for (j in 1:ymax) {
      count[i,j]<-length(which(simul.data.2[,1]>=i-1 & 
      simul.data.2[,1]<i & simul.data.2[,2]>=j-1 & simul.data.2[,2]<j))
  }
}
sum(count)
  
network<-array(0,c(xmax,ymax,n))
edge<-array(0,c(xmax,ymax,n))
y<-rep(NA,n)
m<-rep(NA,n)
newnet<-array(0,c(xmax,ymax,n))
look<-matrix(0,4,2)
  
for (i in 1:n) {
newnet[rand.ind[i,1],rand.ind[i,2],i]<-1
  while (!identical(newnet[,,i], network[,,i])) {
    network[,,i] <- newnet[,,i]
    net.ind<-arrayInd(which(network[,,i]==1), dim(network[,,i]))
      if (count[rand.ind[i,1],rand.ind[i,2]]>0) {
        for (j in 1:length(which(network[,,i]==1))) {
          look[1,]<-c(net.ind[j,1]-1,net.ind[j,2])
          look[2,]<-c(net.ind[j,1]+1,net.ind[j,2])
          look[3,]<-c(net.ind[j,1],net.ind[j,2]-1)
          look[4,]<-c(net.ind[j,1],net.ind[j,2]+1)
          for (k in 1:4) {
            if (0<look[k,1] & look[k,1]<=20 & 0<look[k,2] & look[k,2]<=20) {
              if (count[look[k,1],look[k,2]]==0) {edge[look[k,1],look[k,2],i]<-1} 
            } 
            if (0<look[k,1] & look[k,1]<=20 & 0<look[k,2] & look[k,2]<=20) {
              if(count[look[k,1],look[k,2]]>0) {newnet[look[k,1],look[k,2],i]<-1} 
            }
            
          }
        }
      }
    }
    y[i]<-sum(count[which(network[,,i]==1)])
    m[i]<-sum(network[,,i])
  }
  
  w<-y/m
  mu.acs[h]<-sum(w)/n
  var.acs[h]<-(1-n/N)*sum((w-mu.acs)^2)/(n*(n-1))
  y.srs<-count[rand.ind]
  mu.srs[h]<-sum(y.srs)/n
  var.srs[h]<-(1-n/N)*sum((y.srs-mu.srs)^2)/(n*(n-1))
}

Mu.acs<-N*mu.acs
Var.acs<-N^2*var.acs
Mu.srs<-N*mu.srs
Var.srs<-N^2*var.srs
Mu.acs.avg<-sum(Mu.acs)/100
Mu.acs.spread<-sum((Mu.acs-Mu.acs.avg)^2)/100
Var.acs.avg<-sum(Var.acs)/100
Var.acs.spread<-sum((Var.acs-Var.acs.avg)^2)/100
Mu.srs.avg<-sum(Mu.srs)/100
Mu.srs.spread<-sum((Mu.srs-Mu.srs.avg)^2)/100
Var.srs.avg<-sum(Var.srs)/100
Var.srs.spread<-sum((Var.srs-Var.srs.avg)^2)/100
plot(cluster.1, xlim=c(0,20), ylim=c(0,20), xlab="Spatial clusters (x)",
ylab="Spatial clusters (y)", main="Simulated populations(Biv. Normal, SD=1)")
points(cluster.2, xlim=c(0,20), ylim=c(0,20), col=2)
points(cluster.3, xlim=c(0,20), ylim=c(0,20), col=3)
points(cluster.4, xlim=c(0,20), ylim=c(0,20), col=4)
points(cluster.5, xlim=c(0,20), ylim=c(0,20), col=5)
sum(frame[,1])
\end{verbatim}

\end{document}